\documentclass[a4paper, 12pt]{article}

\usepackage[utf8]{inputenc}
\usepackage[T1]{fontenc}
\usepackage{textcomp}
\usepackage{graphicx}
\usepackage{amsfonts}
\usepackage{enumerate}
\usepackage{amsmath}
\usepackage{amssymb}

\def\C{{\mathbb C}}

\def\Q{{\mathbb Q}}

\def\cG{{\cal G}}
\def\cK{{\cal K}}

\def\cP{{\cal P}}

\def\cK{{\cal K}}

\title{The interlace polynomial of binary delta-matroids and link invariants}

\author{Nadezhda Kodaneva
	\thanks{National Research University Higher School of Economics, Moscow, Russia, e-mail:
		n.kodaneva@ya.ru}}

\begin{document}

	\maketitle

	\begin{abstract}
		
		In this work, we study the interlace polynomial as a generalization of a graph invariant to delta-matroids. We prove that the interlace polynomial satisfies the four-term relation for delta-matroids and determines thus a finite type invariant of links in the $3$-sphere.
		
	\end{abstract}

	
	Finite type invariants of knots were introduced by V.~Vassiliev~\cite{V90} in the end of 1980'ies.
	They can be described in terms of weight systems, which are functions defined on chord diagrams and satisfying the so-called
	four-term relations.
	A graph can be associated to every chord diagram, called the intersection graph,
	and this association leads to introducing the four-term relation for abstract graphs, see~\cite{L00}.
	Each graph invariant satisfying $4$-term relations defines therefore a weight system.
	
	A chord diagram can be interpreted as an embedded graph with a single vertex.
	Considering arbitrary ribbon graphs, not necessarily those having a single vertex, we obtain a generalization of the notion of weight systems from knots to links. 
	To arbitrary embedded graph, no intersection graph can be associated. Instead, one associates to an embedded graph
	another combinatorial structure, namely, its delta-matroid~\cite{B87}, and this delta-matroid is binary. 
	As it was shown in~\cite{LZ17}, it is possible to define the four-term relation for binary delta-matroids 
	in such a way that the mapping taking an embedded graph to its
	delta-matroid respects the corresponding four-term relations. This result raises the question whether known delta-matroid invariants
	respect the introduced relation. In the present paper, we answer in affirmative to this question for the interlace polynomial.
	
	The article was prepared within the framework of the Academic Fund Program at the National Research University Higher School of Economics (HSE University) in 2020 — 2021 (grant №20-04-010) and within the framework of the Russian Academic Excellence Project «5-100».
	
	\section{Introduction}
	
	In this section we introduce the main notions of our study and explain the result of N.~Netrusova 
	for the interlace polynomial of graphs.
	
	\subsection{The interlace polynomial for graphs}
	
	The interlace polynomial was originally introduced as a function defined recursively on two-in two-out directed graphs. It appeared in~\cite{ABS00}, a work on DNA sequencing written by R. Arratia, B. Bollobas, and G.B. Sorkin, and then the recursion was generalized to determine the polynomial of an arbitrary simple graph. Here we will reproduce the definition of the interlace polynomial for graphs following the terminology of \cite{M16}. We first require the definition of the pivot operation.
	
	Let $G$ be a simple graph (that is, a graph without loops and multiple edges). For any pair of adjacent vertices $a, b$ of a graph $G$, the remaining vertices of the graph can be split into four classes:
	
	\begin{enumerate}
		\item vertices adjacent to $a$ and not to $b$;
		\item vertices adjacent to $b$ and not to $a$;
		\item vertices adjacent to both $a$ and $b$;
		\item vertices adjacent neither to $a$ nor to $b$.
	\end{enumerate}
	
	Then the {\it pivot} $G^{ab}$ is the graph obtained from $G$ by erasing existing edges between the vertices in the first three classes if and only if they belong to different classes and adding such edges if they are not present in $G$.

	{\bf Definition 1.} Let $G$ be a graph without loops and multiple edges. {\it The interlace polynomial} $q(G, x)$ is a polynomial in one variable $x$ satisfying the following formulas.
	
	\begin{enumerate}
		\item If the graph $G$ does not have edges, then $q(G, x)=x^n$, where $n$ is the number of vertices of $G$.
		\item If $ab$ is an edge in $G$, then $$q(G, x)=q(G \setminus a,x)+q(G^{ab} \setminus b, x),$$
		where $G \setminus a$ denotes the graph obtained from $G$ by deleting the vertex $a$ and all the edges connecting $a$ with the other vertices.
	\end{enumerate}
	
	A theorem in~\cite{ABS00} states that the interlace polynomial is well-defined: it does not depend on the order in which the recurrence is applied.
	It follows from the definition that $q(G,x)$ is a polynomial of degree $n=|V(G)|$.
	
	
	\subsection{Delta-matroids}
	
	In this section, we give general definitions related to delta-matroids and introduce the way the interlace polynomial can be generalized to delta-matroids. Our exposition follows~\cite{LZ17}.

	{\bf Definition 2.} {\it A set system} $(E, \Phi)$ is a pair consisting of a finite set $E$ and a subset $\Phi$ of $2^E$.
	The set $E$ is the {\it ground set} of the set system. The elements of $\Phi$ are the {\it feasible sets}. The set system is called {\it proper} if $\Phi$ is nonempty.
	
	In this work, we consider only proper set systems.
	
	{\bf Definition 3.} A {\it delta-matroid} is a proper set system satisfying the following {\it symmetric exchange axiom}:  for all $X, Y \in \Phi$, if $x \in X \Delta Y$, then there exists $y \in X \Delta Y$ such that $\{x, y\} \Delta X \in \Phi$. Here $\Delta$ denotes the set symmetric difference, $X\Delta Y=(X\setminus Y)\cup (Y\setminus X)$.
	
	For example, the set system $(\{1, 2, 3\}, \{\emptyset, \{1\}, \{1, 2\}, \{2, 3\}, \{1, 2, 3\}\})$ is a delta-matroid. But the set system $(\{1, 2, 3\}, \{\emptyset,  \{1, 2\}, \{2, 3\}, \{1, 2, 3\}\})$  is not a delta-matroid.
	
	{\bf Definition 4.} Let $D=(E, \Phi)$ be a set system. For $X \subset E$, define the {\it twist} $D * X$  as the set system $(E, \Phi \Delta X)$, where $\Phi \Delta X = \{F \Delta X | F \in \Phi \}$.
	
	According to A.~Bouch\'et, any twist of a delta-matroid is a delta-matroid.
	
	For an abstract graph $G$, we can define its {\it nondegeneracy delta-matroid} $M_{G}$ in the following way. Let $A_G$ denote the adjacency matrix of $G$ over the field $\mathbb{F}_2$. We say that~$G$ is {\it nondegenerate\/} if the matrix $A_G$ is nondegenerate.
	The empty graph is nondegenerate by convention.
	By definition, the ground set of $M_G$ is the set of the graph's vertices $V(G)$, and the
	feasible sets are the subsets $X \subseteq V(G)$ such that the induced subgraph $G|_X$ is nondegenerate.
	
	It can be shown that set systems defined in this way are indeed delta-matroids. Furthermore, we can analogously define a set system $M_{A}$ for an arbitrary symmetric matrix $A$ over the field $\mathbb{F}_2$ (meaning that we allow ones, not only zeroes, on the matrix diagonal). This set system is a delta-matroid as Bouchet showed in~\cite{B87}.  Those delta-matriods that can be obtained as twists
	of delta-matroids of the form $M_A$ are called {\it binary} delta-matroids.
	
	{\bf Definition 5.} An element $e\in E$ of the ground set~$E$ of a delta-matroid $D=(E, \Phi)$ is called {\it a loop} if it is not contained in any feasible set and {\it a coloop} if it is contained in all feasible sets. For $e$ that is not a loop, we define {\it $D$ contract $e$} as the set system $D/e=(E\setminus\{e\}, \{F\setminus \{e\}| F \in \Phi, e \in F\})$. For $e$ that is not a coloop, we define {\it $D$ delete $e$} as the set system $D\setminus e = (E\setminus \{e\}, \{ F \in \Phi | e \notin F \}).$ If $e$ is a loop, then, by definition, $D/e=D\setminus e.$ If $e$ is a coloop, then $D\setminus e=D/e.$
	
	These operations are well-defined on delta-matroids, i.e. $D$ delete $e$ and $D$ contract $e$ are delta-matoids.
	
	{\bf Definition 6.} \cite{BH14} The {\it distance} from a set system $D=(E, \Phi)$ to a subset $X \subseteq E$ is
	$$d_{D}(X)= \min_{F\in\Phi}{|F \Delta X|}.$$
	The {\it interlace polynomial} of $D$ is $$q_{\Delta}(D, x) = \sum\limits_{\phi\subseteq E}x^{d_{D}(\phi)}.$$
	
	The interlace polynomial for the class of delta-matroids associated to graphs coincides with the interlace polynomial for graphs. The following result describes the connection between the polynomials.
	
	{\bf Theorem 1.} \cite{BH14} Let $G$ be a graph, and $M_{G}$ be the nondegeneracy delta-matroid of $G$. Then $q(G,x)=q_{\Delta}(M_{G}, x-1)$.
	
	The recursive formula defining the interlace polynomial for graphs can be also generalized to the case of delta-matroids.
	
	\textbf{Theorem 2.} \cite{BH14} Let $D=(E, \Phi)$ be a delta-matroid with an element $e\in E$ that is neither a coloop nor a loop. Then we have $$q_{\Delta}(D, x)=q_{\Delta}(D\setminus e, x)+q_{\Delta}(D*e \setminus e, x).$$
	If $\emptyset \in \Phi$, then for any $X \subseteq E$ and for any $e \in X$  the equality $$q_{\Delta}(D, x)=q_{\Delta}(D\setminus e, x)+q_{\Delta}(D*X \setminus e, x)$$ holds.
	
	\subsection{The 4-term relation for graphs}
	
	The connection between the interlace polynomial for graphs and the theory of knots was studied by N. Netrusova \cite{Nu}. She showed that the interlace polynomial for graphs satisfies the four-term relation and therefore determines a knot invariant.
	
	{\bf Definition 7.} The algebra $\mathcal{G}$ of graphs is the infinite dimensional commutative algebra (over the field $\mathbb{C}$) freely spanned by all graphs with the multiplication defined as the disjoint union of graphs.
	
	This algebra is graded, i.e. $\mathcal{G}=\mathcal{G}_0 \oplus \mathcal{G}_1 \oplus \dots$, where $\mathcal{G}_{k}$ is the subspace of the algebra spanned by graphs with $k$ vertices.
	
	With this construction, we can consider the interlace polynomial as a function on the algebra $\mathcal{G}$ extending it by linearity. It will be defined correctly since the value of the interlace polynomial on graphs' product is equal to the product of the polynomials on the factors as it is stated in the following proposition.
	
	{\bf Proposition 1.} If a graph $G$ can be represented as a disjoint union of two graphs $G_{1}$ and $G_{2}$, then $q(G, x)=q(G_{1}, x) \cdot q(G_{2}, x).$
	
	{\bf Definition 8.} A function $f$ defined on graphs is called {\it a 4-invariant} if it satisfies the following property for any graph $G$ and any pair of its vertices $a, b \in V(G)$: $$f(G)-f(G'_{ab})= f(\widetilde{G_{ab}}) - f(\widetilde{G'_{ab}});$$
	here $G'_{ab}$ denotes the graph obtained from $G$ by switching the adjacency of $a$ and $b$, $\widetilde{G_{ab}}$ is the graph obtained by switching the adjacency to $a$ of every vertex $c\in V(G)$ if and only if $c$ is connected to $b$. Note that these operations commute and their composition
	produces the graph $\widetilde{G'_{ab}}$. This relation was introduced by S.~K.~Lando in~\cite{L00}.
	
	{\bf Theorem 3.} \cite{Nu} The interlace polynomial of graphs is a 4-invariant.
	
	This result brings us to the question whether the generalization of the interlace polynomial to binary delta-matroids possesses an analogous property.

	\section{The $4$-term relation for the interlace polynomial for delta-matroids}
	
	In this section we recall the $4$-term relation for binary delta-matroids introduced in~\cite{LZ17}
	and prove that the interlace polynomial of delta-matroids satisfies this relation.
	
	\subsection{The $4$-term relation for delta-matroids}
	
	The way the four-term relation is introduced for delta-matroids is motivated by the relation for ribbon graphs. 
	Namely, to an embedded graph~$\Gamma$ its delta-matroid $M_\Gamma$ can be associated in the following way,
	due to A.~Bouch\'et~\cite{B87}. The ground set of~$M_\Gamma$ is the set~$E(\Gamma)$ of edges of~$\Gamma$, and a subset
	$X\subset E(\Gamma)$ is feasible if the embedded induced subgraph $\Gamma|_X$ has a connected boundary (= a unique face).
	
	\textbf{Remark.} In~\cite{KS16}, the $4$-term relations for Lagrangian subspaces were introduced, which also can
	be considered as a combinatorial counterpart of the $4$-term relations for embedded graphs. V.~Zhukov~\cite{Z17}
	proved the equivalence of this approach to the one for binary delta-matroids.
	
	To give the description of four-term relation for $\Delta$-matroids, we should define two Vassiliev moves. The second Vassiliev move (handle sliding) was introduced and studied by I. Moffatt and E. Mphako-Banda in \cite{MMB17}. The first Vassliliev move (exchanging of handle ends) was defined by S. Lando and V. Zhukov in \cite{LZ17}.
	
	\textbf{Definition 9.} Let $D=(E,\Phi)$ be a set system, $a$ and $b$ be two distinct elements in $E$. Then the result of \textit{handle slide taking a over b} is the set system $\widetilde{D_{ab}}=(E, \widetilde{\Phi_{ab}}),$ where $\widetilde{\Phi_{ab}}=\Phi \Delta \{X \sqcup {a} | X \sqcup {b} \in \Phi \textrm{ and } X \subseteq E \setminus \{a, b\}\}$.
	
	\textbf{Proposition 2.} \cite{MMB17} The set of binary delta-matroids is closed under the operation of handle sliding.
	
	\textbf{Definition 10.} Let $D=(E,\Phi)$ be a set system, $a$ and $b$ be two distinct elements in $E$. The result of \textit{exchanging of the ends of the handles $a$ and $b$} is the set system $D'_{ab}= (E, \Phi'_{ab}),$ where $\Phi'_{ab}=\Phi \Delta \{X \sqcup \{a, b\} | X \in \Phi \textrm{ and } X \subseteq E \setminus \{a, b\}\}$.

	\textbf{Proposition 3.} \cite{LZ17} The set of binary delta-matroids is closed under the operation of exchanging the handle ends.
	
	Then a function $f$ defined on binary $\Delta$-matroids is said to satisfy the \textit{four-term relation} if the equality $$ f(D)-f(D'_{ab})= f(\widetilde{D_{ab}}) - f(\widetilde{D'_{ab}})$$ holds for every delta-matroid $D=(E, \Phi)$ and every pair $a$, $b \in E,$ where $a \neq b$.

	\subsection{Interlace polynomial and the $4$-term relation}
	
	Our first main result is the following
	
	\textbf{Theorem 4.} The interlace polynomial for binary delta-matroids satisfies the 4-term relation.
	
	\textit{Proof.}
	Let $D=(E,\Phi)$ be a binary delta-matroid. We claim that for every set $\phi \subseteq E$ one has
	
	\begin{equation}\label{ip4t}
	x^{d_{D}(\phi)}-x^{d_{D'_{ab}}(\phi)}=x^{d_{\widetilde{D}_{ab}}(\phi)}-x^{d_{\widetilde{D}'_{ab}}(\phi)}.
	\end{equation}
	
	This can be proved by showing that for every $\phi$ one of the following conditions holds.
	
	\begin{enumerate}
		\item $d_{D}(\phi)=d_{D'_{ab}}(\phi)=d_{\widetilde{D}_{ab}}(\phi)=d_{\widetilde{D}'_{ab}}(\phi)$
		\item $d_{D}(\phi)=d_{D'_{ab}}(\phi)$ and $d_{\widetilde{D}_{ab}}(\phi)=d_{\widetilde{D}'_{ab}}(\phi)$
		\item $d_{D}(\phi)=d_{\widetilde{D}_{ab}}(\phi)$ and $d_{D'_{ab}}(\phi)=d_{\widetilde{D}'_{ab}}(\phi)$
	\end{enumerate}

	There are four cases:
	
	\begin{enumerate}[I.]
		\item $\phi$ is a feasible set in all four delta-matroids;
		\item $\phi$ is a feasible set only in $D$ and $D'_{ab}$ (or in $\widetilde{D}_{ab}$ and $\widetilde{D}'_{ab}$);
		\item $\phi$ is a feasible set only in $D$ and $\widetilde{D}_{ab}$ (or in $D'_{ab}$ and $\widetilde{D}'_{ab}$);
		\item  $\phi$ is a feasible set in none of the four delta-matroids.
	\end{enumerate}
	
	Other cases are impossible because a set $\phi$ can be feasible in one of the delta-matroids and not be feasible in another one if either $\phi=\phi' \sqcup \{a, b\}$ or $\phi = \phi' \sqcup \{a\},$ where $\phi' \subseteq E \setminus \{a, b\}.$ The first equality corresponds to the first Vassiliev move, and the second one to the second Vassiliev move. Therefore, $\phi$ can be added to the set of feasible sets or removed from it only by one of the moves. And this cannot happen only in one case of applying a move because sets $\phi' \sqcup \{b\}$ and $\phi'$ cannot be added or removed by any of the moves.
	
	Now let us consider all the possibilities case by case.
	
	{\bf I.} All the distances are equal to~$0$, so the equality (\ref{ip4t}) holds.
	
	{\bf II.} If $\phi$ is feasible in $D$ and not feasible in $\widetilde{D}_{ab}$, then $\phi=\phi' \sqcup \{a\}, b\notin \phi$ and $\phi' \sqcup \{b\}$ is feasible in all four delta-matroids. The distances from $\phi$ to $\widetilde{D}_{ab}$ and $\widetilde{D}'_{ab}$ are not greater than $2$, since $|\phi \Delta (\phi' \sqcup \{b\})|=2$. If one of them equals 1, then there exists a set $\phi_{0}$, feasible in $\widetilde{D}_{ab}$ or $\widetilde{D}'_{ab}$, s.t. $\phi \Delta \phi_{0} = c$, $c \in E$. If it is feasible in the other delta-matroid, (\ref{ip4t}) follows. Otherwise, $\phi_{0}=\phi'_{0} \sqcup \{a, b\}$ and $\phi'_{0}$ is feasible in all delta-matroids. Since $b$ does not belong to $\phi$, $c=b$ and $\phi'=\phi'_{0}$. Therefore, $d_{\widetilde{D}_{ab}}(\phi)=d_{\widetilde{D}'_{ab}}(\phi)=|\phi \Delta \phi'|=1.$
	
	{\bf III.} If $\phi$ is feasible in $D$ but not in $D'_{ab}$, then $\phi=\phi' \sqcup \{a, b\}$ and $\phi'$ is feasible in all four delta-matroids. Similarly, the distances from $\phi$ to $D'_{ab}$ and $\widetilde{D}'_{ab}$ are not greater than 2. If one of the distances is equal to 1, then for a set $\phi_0$ that is feasible in $D'_{ab}$ or $\widetilde{D}'_{ab}$, $\phi \Delta \phi_0 = c$, where $c$ is an element of the ground set. If $\phi_0$ is feasible in both delta-matroids, then the other distance also equals 1. Conversely, $\phi_0=\phi_0' \sqcup \{a\}$, $b\notin \phi_0$, and $\phi_0' \sqcup \{b\}$ is feasible in both delta-matroids. In addition, $b=c$, $\phi'=\phi_0'$. Hence, the distance from $\phi$ to $\phi_0$ is the same as the distance from $\phi$ to $\phi_0' \sqcup \{b\}$, i.e. the distances from $\phi$ to both delta-matroids equal 1.
	
	{\bf IV.} Without loss of generality, we can assume that the distance from $\phi$ to $D$ is not greater than the three other distances. Let $\phi_{0}$ denote a feasible set in $\Phi$ such that $d_{D}(\phi)=|\phi \Delta \phi_{0}|$.
	If not all the distances from $\phi$ to the four delta-matroids have the same values, then $\phi_{0}$ is not feasible either in $D'_{ab}$ or in $\widetilde{D}_{ab}$.
	
	In the first case, $\phi_0=\phi'_0 \sqcup \{a, b\}$ and $\phi'_0$ is feasible in all four delta-matroids. If exactly one element out of
	$\{a,b\}$ belongs to $\phi$, then the distances from $\phi$ to $\phi_0$ and $\phi'_0$ are equal. Therefore, the distances from $\phi$ to all four delta-matroids are equal to $|\phi \Delta \phi'_0|$. If neither $a$ nor $b$ belongs to $\phi$, then $|\phi \Delta \phi'_0|<|\phi \Delta \phi_0|$, whence this case is impossible.
	
	If both $a$ and $b$ are in $\phi$, then $|\phi \Delta \phi'_0|=|\phi \Delta \phi_0|+2.$ Suppose that at least one of the distances from $\phi$ to $\widetilde{D}_{ab}$ and $\widetilde{D}'_{ab}$ is less than this value. Then this distance is achieved for some $\phi_1$ either from $\widetilde{\Phi}_{ab}$ or $\widetilde{\Phi}'_{ab}$. If the set $\phi_1$ is feasible in both delta-matroids, then the equality (\ref{ip4t}) holds. Otherwise, $\phi_1=\phi'_1 \sqcup \{a\}$, $b\notin \phi'_1$, and $\phi'_1 \sqcup \{b\}$ is feasible in both delta-matroids. By assumption, $a, b \in \phi$, hence the distances from $\phi$ to $\phi_1$ and $\phi'_1 \sqcup \{b\}$ are equal.
	
	In the second case, $\phi_0=\phi'_0 \sqcup \{a\}$ and $\phi'_0 \sqcup \{b\}$ is feasible in all four delta-matroids. If both $a$ and $b$ belong to $\phi$ or do not belong to $\phi$, we have $|\phi \Delta \phi_0|=|\phi \Delta (\phi_0' \sqcup \{b\})|$. Therefore, the distance from $\phi$ to $D$ equals three other distances. If $b$ belongs to $\phi$ and $a$ does not belong to $\phi,$ then $|\phi \Delta (\phi_0' \sqcup \{b\})|=|\phi \Delta \phi_0|-2$ what is impossible under the assumptions. Otherwise, $\phi = \phi' \sqcup \{a\}$ and $b \notin \phi$. Then the distances from $\phi$ to $D'_{ab}$ and $\widetilde{D}'_{ab}$ are not greater than $|\phi \Delta (\phi_0' \sqcup \{b\})|=|\phi \Delta \phi_0|+2.$  If one of the distances is less than this value, it is reached for a feasible set $\phi_1$ and the other distance has the same value. Indeed, either $\phi_1$ is feasible in the other delta-matroid or $\phi_1=\phi_1' \sqcup \{a, b\}$, $\phi_1'$ is also feasible, and the number of elements in $\phi \Delta \phi_1$ is the same as in $\phi \Delta \phi_1'$.

	Thus, for every set $\phi$ we have either four equal distances or two pairs of equal distances.
	
	\hfill $\square$
	
	It was shown in [5] that any invariant of binary delta-matroids satisfying the $4$-term relations
	determines a finite type invariant of links. As a consequence, we obtain the following
	
	\textbf{Corollary.} The interlace polynomial of binary delta-matroids defines a finite type invariant of links in the 3-sphere.

	\section{The interlace polynomial for primitive elements in the Hopf algebra of delta-matroids}
	
	In this section we study the distinguishing power of the interlace polynomial of binary delta-matroids.
	
	Binary delta-matroids span a commutative cocommutative Hopf algebra. This Hopf algebra is an algebra
	of polynomials in its primitive elements, and
	it makes sense to consider the restriction of a multiplicative invariant
	to the subspace of primitives (which determines the
	polynomial completely). In particular, the distinguishing power of a multiplicative invariant
	can be measured as the dimension of the space of its values on 
	the subspace of primitive elements in degree~$n$.
	For the interlace polynomial for graphs, N.~Netrusova~\cite{Nu} proved
	that its space of values on primitive elements in degree~$n$ is $\left[\frac{n}2\right]$,
	for $n\ge2$. (This statement can be compared with the similar statement for the
	chromatic polynomial, for which the dimension in question is~$1$, for all $n=1,2,3,\dots$.)
	
	Our main result here is that the dimension of the space of values of the interlace polynomial
	for binary delta-matroids in degree~$n$ is $n$, for all~$n$.
	
	\subsection{The Hopf algebra of binary delta-matroids}
	
	Let us define the product and the coproduct of binary delta-matroids on the vector space $\mathcal{B}$ freely spanned by all binary delta-matroids. The space $\mathcal{B}$ is graded, i. e. $\mathcal{B}=\mathcal{B}_{0} \oplus \mathcal{B}_{1} \oplus \mathcal{B}_{2}  \oplus \dots$ where $\mathcal{B}_{k}$ is the subspace spanned by binary delta-matroids with $k$ elements in the ground set.
	
	The \textit{product} of delta-matroids $D_{1}=(E_{1}, \Phi_{1})$ and $D_{2}=(E_2, \Phi_2)$ is the delta-matroid $D=(E_1 \sqcup E_2, \Phi)$ where $\Phi = \{\phi_1 \sqcup \phi_2 | \phi_1 \in \Phi_1 \textrm{ and } \phi_2 \in \Phi_2\}$.
	
	This multiplcation extended by linearity to the space $\mathcal{B}$ makes it into a commutative algebra. Its unit is the delta-matroid with an empty ground set. The algebra is graded since the multiplication $m$ preserves the grading: $$m:\mathcal{B}_k \otimes \mathcal{B}_l \rightarrow \mathcal{B}_{k+l}. $$
	The \textit{coproduct} of a delta-matroid $D=(E, \Phi)$ is 
	$$\mu(D)=\sum\limits_{E'\subseteq E} D_{E'} \otimes D_{E\setminus E'}, $$ 
	where $D_{E'}=D\setminus(E\setminus E')$ is the \textit{restriction of D to E'}.
	

	\textbf{Definition 11.} An element $p$ of a bialgebra is said to be \textit{primitive}, if $\mu(p)=1\otimes p+p\otimes 1$.
	
	Primitive elements form a subspace of $\mathcal{B}$. Each homogeneous subspace~$\mathcal{B}_n$ is a direct sum
	of its primitive subspace and the subspace of decomposable elements, which is spanned by the products
	of elements of smaller degree. This decomposition determines a projection of~$\mathcal{B}_n$ to its subspace
	of primitive elements. The projection of a delta-matroid $D=(E, \Phi)$ to the subspace of primitive elements
	is as follows,
	see~\cite{L99}:
	$$\pi(D)=D-1! \sum\limits_{E_1 \sqcup E_2 = E} D_{E_1} D_{E_2} +2! \sum\limits_{E_1 \sqcup E_2 \sqcup E_3 = E} D_{E_1} D_{E_2} D_{E_3} -  \dots $$
	where the $k^{\textrm{th}}$ sum is taken over all partitions of the ground set $E$ into $k$ nonempty subsets $E_1$, \dots ~$E_k$.
	
	

	\subsection{Interlace polynomial for primitives}
	
	We say that a polynomial selects $k$ primitive elements in degree $n$, if it has $k$ linearly independent values on the subspace of primitive elements with $n$ elements in the ground set.
	Our second main result is that the interlace polynomial selects $n$ primitive elements in degree $n$. The following lemmas are needed for the sequel.

	\textbf{Lemma 1.} Let $D=(E, \Phi)$ be the nondegeneracy delta-matroid associated to an $n\times n$-matrix $A=((a_{ij}))_{i, j=1}^{n}$
	and let $A'$ be the $(n+1)\times(n+1)$-matrix obtained from $A$ by adding the $(n+1)$~th row and
	the $(n+1)$~th column to the~matrix~$A$, with $a_{1, n+1}=a_{n+1, 1}=1$ and $a_{i, n+1}=a_{n+1, i}=0$ for all $i$ from $2$ to~$n+1$.
	Let $D'$ be the~nondegeneracy delta-matroid of~$A'$. Then $q(D')=q(D)+(x+1)q(D\setminus 1)$.
	
	\textit{Proof.}
	It follows from Theorem 2 that $q(D')=q(D'\setminus \{n+1\})+q((D' * \{n+1\}) \setminus \{n+1\})$.
	Clearly, $D'\setminus \{n+1\} = D$. It remains to check that $q((D' * \{n+1\}) \setminus \{n+1\}) = (x+1) q(D\setminus 1)$.
	
	
	A set $\phi$ is feasible in $(D' * \{n+1\}) \setminus \{n+1\}$ if and only if $\phi$ does not contain $n+1$ and $\phi \sqcup \{n+1\}$ is feasible in $D'$. Every such set contains $1$ since a principal submatrix of $A'$ including the last row and  not including the first column has only zeros in the last row and the set corresponding to such submatrix cannot be feasible.
	
	Moreover, a set $\phi \sqcup \{1\}$ is feasible in $(D' * \{n+1\}) \setminus \{n+1\}$ if and only if the set $\phi$ is feasible in $D\setminus 1$. Indeed,
	$$ \det \begin{pmatrix}
	a_{11} & \dots & 1\\
	\dots & A[\phi] & 0\\
	\dots & & \dots\\
	1 & 0\dots & 0\\
	\end{pmatrix}
	\quad = \det
	\begin{pmatrix}
	a_{11} \dots & 1\\
	& 0\\
	A[\phi] & \dots\\
	&  0\\
	\end{pmatrix} =
	\det (A[\phi]).
	$$
	
	Then 
	\begin{eqnarray*}
		q((D' * \{n+1\}) \setminus \{n+1\}) &=&\sum\limits_{X \subseteq E} x^{d_{(D' * \{n+1\}) \setminus \{n+1\}} (X)} \\
		&=&  \sum\limits_{X \subseteq E\setminus \{1\}} (x^{d_{(D' * \{n+1\}) \setminus \{n+1\}} (X)} + x^{d_{(D' * \{n+1\}) \setminus \{n+1\}} (X\sqcup \{1\})}) \\
		&=&  \sum\limits_{X \subseteq E\setminus \{1\}} (x^{d_{D \setminus \{1\}} (X)+1} + x^{d_{D \setminus \{1\}} (X)})\\
		&=& \sum\limits_{X \subseteq E\setminus \{1\}} x^{d_{D \setminus \{1\}} (X)} (x+1)\\
		&=& q(D\setminus \{1\}) (x+1).
	\end{eqnarray*}

	\hfill $\square$
	
	Below, we will denote $(D^{[k]})'$ by $D^{[k+1]}$, with $D^{[0]}=D$.
	
	\textbf{Lemma 2.} Let $D_{n}$ be the delta-matroid $(\{1, \dots n\}; \{\emptyset, \{1\}, \dots, \{n\}\})$; then the interlace polynomial of its projection onto the subspace of primitive elements has a non-zero constant term. (It equals $(-1)^{n-1} (n-1)!$.)
	
	\textit{Remark.} I am grateful to S.~K.~Lando for the proof of the fact that the sum that appears
	in the constant term below indeed equals  $(-1)^{n-1} (n-1)!$.

	\textit{Proof.}
	
	The delta-matroid $D_{n}$ is a nondegeneracy delta-matroid. It is associated to the $n\times n$-matrix  with every element equal to one. Consider the projection of $D_n$ to the subspace of primitives,
	$$\pi(D_n)=D_n-\sum\limits_{E_1\sqcup E_2 = E} D_{E_1} D_{E_2} + \dots$$.
	
	Each delta-matroid $D_{E_i}$ also is a nondegeneracy delta-matroid and the corresponding 
	matrix is the $|E_i| \times |E_i|$-matrix where every element is equal to one. Therefore, $D_{E_i} = D_{|E_i|}$.
	
	For an~arbitrary $k$ from $1$ to $n$, the delta-matroid $D_k$ has $k+1$ feasible sets, namely,
	$\emptyset$, $\{1\}$, \dots, $\{k\}$. Therefore, the~constant term of~$q(D_k)$ equals $k+1$. Indeed, a~summand of~$\sum\limits_{X \subseteq \{1,\dots k\}} x^{d_{D_{k}}}(X)$ has x in the~degree $0$ if and only if the~set $X$ is feasible.
	
	Now, the~constant term of~$q(\pi(D_n))$ equals $$
	n+1-1!\sum_{{E_1\sqcup E_2=E\atop E_1\ne\emptyset,E_2\ne\emptyset}}(|E_1|+1)(|E_2|+1)
	+2!\sum_{{E_1\sqcup E_2\sqcup E_3=E\atop E_1\ne\emptyset,E_2\ne\emptyset,E_3\ne\emptyset}}(|E_1|+1)(|E_2|+1)(|E_3|+1)-\dots.
	$$
	
	Let us prove that this sum is equal to $(-1)^{n-1}(n-1)!$ for $n\geqslant 2$.

	Consider the subalgebra~$\cK$ spanned by complete graphs in the Hopf algebra of graphs~$\cG$ over the field~$\Q$. Denote the complete graph on~$n$ vertices by~$K_n$, $n=0,1,2,\dots$ and introduce 
	the exponential generating function
	
	\begin{eqnarray*}
		\cK(t)&=&\sum_{n=0}^\infty K_n\frac{t^n}{n!}\\
		&=&1+K_1\frac{t^1}{1!}+K_2\frac{t^2}{2!}+\dots
	\end{eqnarray*}
	(here we use the fact
	that the~empty graph $K_0$ is the unit in the Hopf algebra of graphs).
	
	The~projection $\pi:\cK\to\cP(\cK)$ to the~subspace of~primitive elements whose kernel is
	the subspace of decomposable elements in~the~Hopf algebra of complete graphs~$\cK$ yields
	\begin{eqnarray*}
		\pi(\cK(t))&=&\sum_{n=0}^\infty \pi(K_n)\frac{t^n}{n!}\\
		&=&\log\cK\\
		&=&K_1\frac{t^1}{1!}+(K_2-K_1^2)\frac{t^2}{2!}+(K_3-3K_1K_2+2K_1^3)\frac{t^3}{3!}+\dots
	\end{eqnarray*}
	
	Here the absolute value of the coefficient of $K_1^{m_1}K_2^{m_2}\dots$ equals 
	the number of partitions of $E$ into~$m_1$ $1$-element subsets, $m_2$ $2$-element subsets, etc. 
	And these are the partitions that appear on the left hand side of the equality.
	
	Consider the~multiplicative mapping  $\varphi:\cK\to\C$  whose value on~the~complete graph~$K_n$ is $n+1$. For~this~mapping,
	\begin{eqnarray*}
		\varphi\cK(t)&=&\sum_{n=0}^\infty\varphi(K_n)\frac{t^n}{n!}\\
		&=&1+\varphi(K_1)\frac{t^1}{1!}+\varphi(K_2)\frac{t^2}{2!}+\dots\\
		&=&1+2\cdot\frac{t^1}{1!}+3\cdot\frac{t^2}{2!}+4\cdot\frac{t^3}{3!}+\dots\\
		&=&\left(\frac{t^1}{0!}+\frac{t^2}{1!}+\frac{t^3}{2!}+\frac{t^4}{3!}+\dots\right)'\\
		&=&(te^t)'\\
		&=&(1+t)e^t.
	\end{eqnarray*}
	
	For the~logarithm of the~generating function~$\varphi\cK$, we have
	\begin{eqnarray*}
		\log\varphi\cK(t)&=&\varphi\log\cK(t)\\
		&=&\sum_{n=0}^\infty\varphi(\pi(K_n))\frac{t^n}{n!}
	\end{eqnarray*}
	(the~first equality holds since the~map~$\varphi$ is multiplicative).
	
	On the other hand,
	\begin{eqnarray*}
		\log\varphi\cK(t)&=&\log((1+t)e^t)\\
		&=&t+\log(1+t)\\
		&=&t+\sum_{n=1}^\infty(-1)^{n-1}\frac{t^n}{n}\\
		&=&t+\sum_{n=1}^\infty(-1)^{n-1}(n-1)!\frac{t^n}{n!}.
	\end{eqnarray*}
	
	Comparing these two formulas of~the~logarithm of the~generating function~$\varphi\cK$, we obtain that the equality holds and $\pi(D_n)$ has a non-zero constant term.
	
	\hfill $\square$

	\textbf{Lemma 3.} If $D$ is a nondegeneracy delta-matroid and $D'$ is obtained from $D$ as in Lemma 1, then $q(\pi (D'))=-x\cdot q(\pi (D))$.
	
	\textit{Remark.} The proof of this lemma is analogous to the proof of similar statement for graphs in [8].
	
	\textit{Proof.}
	
	Partitions of $E\sqcup \{n+1\}$ into nonempty subsets can be divided into three kinds: the ones having $n+1$ as a separate subset,  the ones where $1$ and $n+1$ are in the same subset, and the ones where $1$ and $n+1$ are in different subsets and $n+1$ is not a separate subset. Assume that $1$ belongs to $E_1$.
	
	In the first case, the delta-matroid corresponding to a partition $E_1\sqcup~E_2\sqcup\dots \sqcup E_k$ is 
	$(D')_{E_1} (D')_{E_2} \dots (D')_{E_k}$. It is isomorphic to the delta-matroid $(D')_{E_1}D_{E_2}\dots D_{E_k}$ since $n+1$ and $1$ belong to $E_1$. Now, the interlace polynomial of such delta-matroid is
	\begin{eqnarray*}
		q((D')_{E_1}D_{E_2}\dots D_{E_k})
		&=&
		q((D')_{E_1}) \cdot q(D_{E_2}) \dots q(D_{E_k})\\
		&=&
		(q(D_{E_1}) + (x+1)q(D_{E_1} \setminus 1)) \cdot q(D_{E_2}) \dots q(D_{E_k})\\
		&=&
		q(D_{E_1}\dots D_{E_k}) + (x+1) \cdot q(D_{E_1 \setminus 1} D_{E_2} \dots D_{E_k}).
	\end{eqnarray*}
	
	By $D_{\{n+1\}}$, we denote the restriction of $D'$ to the set $\{n+1\}$ or, equivalently, the delta-matroid $(\{n+1\}, \{\emptyset\})$. For an arbitrary delta-matroid $M$ having $n$ elements in the ground set, the distance from a set to the delta-matriod $MD_{\{n+1\}}$ is equal to the distance from this set to the delta-matroid $M$ because the only feasible set in $D_{\{n+1\}}$ is the empty set and a set is feasible in $M$ if and only if it is feasible in $MD_{\{n+1\}}$.
	Now,
	\begin{eqnarray*}
		q(MD_{\{n+1\}})
		&=&
		\sum\limits_{X \subseteq \{1,\dots n+1\}} x^{d_MD_{\{n+1\}}}(X)\\
		&=&
		\sum\limits_{X \subseteq \{1,\dots n\}} (x^{d_M}(X) + x^{d_M}(X \sqcup \{n+1\}))\\
		&=&
		\sum\limits_{X \subseteq \{1,\dots n\}} x^{d_M}(X)(1 + x)\\
		&=&
		(x+1)q(M).
	\end{eqnarray*}
	
	Finally, if $E_i$ is the subset containing $n+1$, then $D_{E_i \sqcup \{n+1\}}=D_{E_i} D_{\{n+1\}}$. Indeed, $n+1$ is not en element of any feasible set in $D_{E_i \sqcup \{n+1\}}$ so the only feasible set in $D_{\{n+1\}}$ is the empty set and the product of $D_{E_{i}}$ and $D_{\{n+1\}}$ has the same feasible sets as $D_{E_{i}}$ and $D_{E_{i} \sqcup \{n+1\}}$. The number of summands $D_{E_1} \dots D_{E_k} D_{\{n+1\}}$ in the projection of $D'$ onto the subspace of primitive elements equals $k-1$ since $n+1$ can be contained in $E_2$, $E_3$ \dots $E_k$.
	
	Now, consider the projection $\pi(D')$ and the value of the interlace polynomial on it,
	\begin{eqnarray*}
		\pi (D')&=& D' - \sum\limits_{E_1\sqcup E_2 = E} D'_{E_1} D_{E_2}
		+
		2! \sum\limits_{E_1\sqcup E_2 \sqcup E_3= E} D'_{E_1} D_{E_2} D_{E_3} -\dots\\
		&&- D D_{\{n+1\}} + 2! \sum\limits_{E_1 \sqcup E_2 =E} D_{E_1} D_{E_2} D_{\{n+1\}} - \dots \\
		&&- \sum\limits_{E_1 \sqcup E_2 =E} D_{E_1} D_{E_2} D_{\{n+1\}} + 2! \cdot 2 \sum\limits_{E_1\sqcup E_2 \sqcup E_3 =E} D_{E_1} D_{E_2} D_{E_3} D_{\{n+1\}} - \dots, 
	\end{eqnarray*}
	and
	\begin{eqnarray*}
		q(\pi(D'))&=&
		q(D - \sum\limits_{E_1\sqcup E_2 = E} D_{E_1} D_{E_2} + \dots)
		+
		(x+1)q(D\setminus 1 - \sum\limits_{E_1\sqcup E_2 = E} D_{E_1 \setminus 1} D_{E_2} \dots) +\\
		&&+ (x+1)q(-D+2! \sum\limits_{E_1\sqcup E_2 = E} D_{E_1} D_{E_2} - \dots) +\\
		&&
		+ (x+1)q(-\sum\limits_{E_1\sqcup E_2 = E} D_{E_1} D_{E_2} +2! \cdot 2 \sum\limits_{E_1\sqcup E_2\sqcup E_3 = E} D_{E_1} D_{E_2} D_{E_3} - \dots) \\ &=&
		q(\pi(D))
		+
		((x+1)q(D\setminus 1 - \sum\limits_{E_1\sqcup E_2 = E} D_{E_1 \setminus 1} D_{E_2} \dots)
		+
		(x+1)(-D + \\
		&& + (2!-1)\sum\limits_{E_1\sqcup E_2 = E} D_{E_1} D_{E_2}
		-
		(3!-2\cdot 2!) \sum\limits_{E_1\sqcup E_2\sqcup E_3 = E} D_{E_1} D_{E_2} D_{E_3} + \dots\\
		&=&
		q(\pi(D))
		+
		((x+1)q(D\setminus 1 - \sum\limits_{E_1\sqcup E_2 = E} D_{E_1 \setminus 1} D_{E_2} \dots)
		+
		(x+1)q(-\pi(D)).
	\end{eqnarray*}
	
	Now, it remains to prove that $$D\setminus 1 - \sum\limits_{E_1\sqcup E_2 = E} D_{E_1 \setminus 1} D_{E_2} +
	2! \sum\limits_{E_1\sqcup E_2\sqcup E_3 = E} D_{E_1 \setminus 1} D_{E_2} D_{E_3} \dots=0.$$

	Every summand except the~first one and the~last one includes partitions of two types. Some of them have at~least one element except $1$ in the~set $E_1$ and the others have $E_1 = \{1\}$. The~first summand consists of one element of~the~first type and the~last one~--- of~the second type.
	If we delete one from a delta-matroid $D_{\{1\} \sqcup E_2 \sqcup \dots E_{k}}$, we will get a delta-matroid $D_{E_2 \sqcup \dots E_{k}}$. There are exactly $k-1$ such delta-matroids of the first type in the previous summand because we can add $1$ to any set of $E_2$, \dots $E_k$ and obtain a partition of $E$ into $k-1$ subsets and a~delta-matriod corresponding to it. In the sum, the~first delta-matriod is multiplied by $(-1)^{k-1} (k-1)!$ and the~others are multiplied by $(-1)^{k-2} (k-2)!$, so their sum equals zero.
	
	Therefore, the sum equals zero and $q(\pi (D'))=-x\cdot q(\pi (D))$.

	\hfill $\square$

	\textbf{Lemma 4.} Let $D_{K_n}$ be the nondegeneracy delta-matroid of the complete graph on $n$ vertices. The degree
	of the polynomial $q(\pi (D_{K_n}))$ equals $n$, and for any $k$ from $0$ to $n-2$ the degree of $q(\pi (D_{n-k}^{[k]}$)) is less than $n$.
	
	\textit{Proof.} The interlace polynomial of an arbitrary delta-matroid is a polynomial with non-negative coefficients. Its power is not greater than $n$ since the number of elements in the symmetric difference of two subsets of the ground set cannot be greater then $|E|=n$.
	
	Note that the empty set is feasible in each of these delta-matroids, whence it is feasible in every summand of their projections. Hence, the distance from a set to a delta-matroid is not greater than the number of elements in the set. It can acquire the value $n$ only for a set with $n$ elements, i.e. for the ground set, and only if there are no feasible sets other than the empty set. Since every two element set is feasible in $D_{K_n}$, the only summand in its projection satisfying this condition is the last summand corresponding to the partition of the ground set into $n$ non-empty subsets. Therefore, the interlace polynomial of the projection $\pi (D_{K_n})$ is a polynomial of degree $n$ and the coefficient of $x^n$ is $(-1)^{n-1}(n-1)!$.
	
	In addition, if a delta-matroid has a one-element feasible set, it remains feasible in every summand of the projection. This implies that the condition for having degree $n$ is not satisfied for any of the delta-matroids $D_{n-k}^{[k]}$. (The set $\{1\}$ is feasible in every of them.)
	
	\hfill $\square$

	\textbf{Theorem 5.} The interlace polynomial for delta-matroids selects $n$ primitive elements in degree $n$.
	
	\textit{Proof.} We claim that the projections of delta-matroids $D_{n}$, $D_{n-1}^{[1]},\dots$ $ D_{2}^{[n-2]} $ and $D_{K_n}$ have linearly independent values of the interlace polynomial.
	
	The proof is by induction on $n$. Assume that $q(\pi(D_{n-1})),~\dots~q(\pi(D_{2}^{[n-3]}))$ are linearly independent. Then by Lemma 3, $q(\pi (D_{n}))=-x \cdot  q(\pi (D_{n-1})), \dots \\ q(\pi (D_{2}^{[n-2]}))=-x \cdot q(\pi (D_{2}^{[n-3]}))$ are linearly independent.  All these polynomials are divisible by $x$, so $q(\pi (D_{n}))$ is linearly independent from them, since it has a non-zero constant term by Lemma 2.
	
	It remains to check that $q(\pi (D_{K_n}))$ is linearly independent from the other polynomials. As Lemma 4 shows, the power of $q(\pi (D_{K_n}))$ is equal to $n$ and other polynomials' power is less than $n$.
	
	\hfill $\square$
	

\end{document}